\documentclass[11pt]{article} 
\usepackage{paralist} 
\usepackage{setspace}
\usepackage{geometry}

\usepackage{amsmath}
\allowdisplaybreaks
\usepackage{amssymb}
\usepackage{amsfonts}
\usepackage{color}
\usepackage{enumerate}

\usepackage{booktabs}

\usepackage{etoolbox}
\patchcmd{\thebibliography}{\section*{\refname}}{}{}{}

\usepackage{verbatim}

\usepackage{graphicx}
\usepackage{subfig}
\usepackage{bm}
\usepackage{tikz}
\bgroup

\usepackage[bookmarks,bookmarksopen=false,
                  pdfauthor={Autor},%
                  pdftitle={Titel},%
                  colorlinks=true,%
                  linkcolor=blue,%
                  citecolor=blue,%
                  urlcolor=blue,%
                ]{hyperref}

\newtheorem{theorem}{Theorem}
\newtheorem{definition}{Definition}
\newtheorem{proposition}[theorem]{Proposition}
\newtheorem{corollary}[theorem]{Corollary}
\newtheorem{conjecture}{Conjecture}
\newtheorem{lemma}[theorem]{Lemma}
\newtheorem{remark}[theorem]{Remark}
\newtheorem{example}{Example}

\newcommand{\proof}{\textbf{Proof:\ }}
\newcommand{\pbox}{\hfill$\Box$\\}

\newcommand{\R}{\mathbb{R}}

\newcommand{\N}{\mathbb{N}}
\newcommand{\C}{\mathbb{C}}
\newcommand{\Z}{\mathbb{Z}}

\renewcommand{\H}{\mathcal{H}}

\renewcommand{\ker}{{\sf Ker}\,}

\newcommand{\ud}{\,\mathrm{d}}

\definecolor{darkviolet}{rgb}{0.58,0,0.83} 

\usetikzlibrary{patterns}

\begin{document}

\title{Sampling trajectories for the short-time Fourier transform  }
\author{  Michael Speckbacher 
\thanks{  Acoustics Research Institute, Austrian
    Academy of Sciences, Wohllebengasse 12-14, 1040 Vienna, Austria, speckbacher@kfs.oeaw.ac.at}}
\date{}
\maketitle


\begin{abstract}
\noindent We study the  problem of stable reconstruction of the short-time Fourier transform from samples taken from trajectories in $\R^2$. 
 We first consider the  interplay between relative density of the trajectory and the reconstruction property. Later, we consider spiraling curves, a special class of trajectories, and connect  sampling and uniqueness properties of these sets. Moreover, we show that for window functions   given by a  linear combination of Hermite functions, it is indeed possible to stably reconstruct from samples on some particular natural choices of spiraling curves.
\end{abstract}

\medskip

\noindent\textbf{MSC2010:} 42C15, 42C30, 42C40, 65T60,  94A20
\newline
\textbf{Keywords:} Gabor analysis, sampling theory, mobile sampling, polyanalytic functions

\section{Introduction}

Reconstructing a function from incomplete data  is an omnipresent task in fields like signal processing or harmonic analysis. Classically, one aims to reconstruct a function from  discrete samples    (see e.g. \cite{seip04,za93})
and the stable version of this problem leads to the notion of frames \cite{christ1,groe1}. 
For a variety of function spaces,   the lower Beurling density \cite{beurl66} (the  definition being adapted to the particular geometry of the function space in consideration) of the discrete sampling set needs to exceed a certain threshold in order for the sampling process to form a frame, see for example  \cite{abba12,la67,rast95,seip93b}. This threshold is commonly known as  the Nyquist rate. 

A more abstract form of the sampling problem is to characterize so called \emph{sampling measures} \cite{Ascensi,lue98,orce98}, i.e. measures $\mu$ on $X$ that satisfy
$$
A\|F\|_{L^p(X)}^p\leq \int_X |F|^p\ud \mu \leq B\|F\|_{L^p(X)}^p,
$$
for some $A,B>0$ and every $F$ in a given subspace of $L^p(X)$. A common situation is  that  $X=\R^d$ and the measure takes the form $\mu=\chi_\Gamma \H^k$, where $\Gamma\subset X$, and $\H^k$, $0\leq k\leq d$, is the $k$-dimensional Hausdorff measure. The case $k=0$ corresponds to the classical discrete sampling problem. For the second extreme case $k=d$, it turns out that a weaker notion of density is  sufficient to stably reconstruct. For the Paley-Wiener space,  this  is known as the Logvinenko-Sereda theorem \cite{kov01,lose73} and  multiple extensions to other spaces of analytic functions
\cite{hakako20,lue81,orce98} and the range space of the short-time Fourier transform \cite{jasp20} have since been established.

Recently, the study of the intermediate cases $0<k<d$ on the Paley-Wiener space $PW_2(\Omega)$ has drawn increased attention and became known in the literature as the mobile sampling problem \cite{unve13,unve13/1}. Mobile sampling has many applications whenever a signal is measured by a  moving sensor like in an MRI scan \cite{bewu00,unve13}. Adapting  Beurling's lower  density for $k=1$, the   \emph{lower path density} of a set $\Gamma$ (measuring the limit average length of the curve $\Gamma$ in  a balls of increasing radii) was studied in 
  \cite{grrounve15}.  In \cite{jami21}, it is shown that a sufficiently large lower path density yields a stable sampling process. On the other hand, there is no  Nyquist rate that the lower path density needs to exceed for stable sampling, see \cite{grrounve15}.
For particularly structured classes of trajectories $\Gamma$ however, Nyquist rates were established in terms of certain parameters that characterize the 'separation'  of $\Gamma$ \cite{janero20,raulzl20}.

In this paper, we study the equivalent of the mobile sampling problem for the \emph{short-time Fourier transform}
$$
V_gf(z)=\int_\R f(t)\overline{g(t-x)}e^{-2\pi i \xi t}\ud t, \quad z=(x,\xi)\in\R^2. 
$$
In particular, if $M^p(\R)$ denotes the modulation space associated to $L^p(\R^2)$ (see e.g. \cite{groe1}), we  study  \emph{Gabor sampling trajectories} for $M^p(\R)$, that is, trajectory sets $\Gamma\subset \R^2$ for which  
\begin{equation}\label{eq:def-intro}
A\|f\|^p_{M^p}\leq \int_\Gamma |V_g f(z)|^p\ud \H^1(z)\leq B\|f\|^p_{M^p},\quad f\in M^p(\R),
\end{equation}
holds. This is a particular  instance of a sampling measure for the short-time Fourier transform \cite{Ascensi}. 

Just like in classical Gabor analysis, the stable sampling property  relies heavily on both the window $g$ and the trajectory set $\Gamma$. On the one hand, one can quickly show that for windows in $M^1(\R)$, a necessary condition for $\Gamma$ being a Gabor sampling trajectory is that it is relatively dense \cite{Ascensi} (see Section~\ref{subsec:nec} for the definition). 
On the other hand, applying a characterization of   sampling measures on the Bargmann-Fock space \cite{orce98} we show that relative density is also sufficient for Gabor sampling trajectories with Gaussian window  (see Section~\ref{subsec:gauss}). It is thus hopeless to search for a  Nyqist rate both in terms of a path density as well as certain separation parameters that holds for general windows. 
 However, in  Proposition~\ref{thm:planar-sampling-frames} we prove that for certain windows 
it is a sufficient for $\Gamma$ being a Gabor sampling trajectory that, for $R$ small enough, $\H^1(\Gamma\cap B_R(z))$ is bounded away from zero.
 
 Later, we change perspective and  study   particular examples of   trajectory sets.  First, we  characterize   sampling and uniqueness properties if $\Gamma$ is a collection of parallel lines. This will in turn  be useful when studying the sampling property of a particular class of trajectories: spiraling curves. Spiraling curves   were introduced in \cite{janero20}. In this article, we work with a slightly more restrictive definition of spiraling curves that nevertheless includes all the examples given in \cite{janero20}. Loosely speaking, a spiraling curve is a trajectory that (in the limit) in  each direction approaches a collection of equispaced parallel lines or equispaced parallel edges. We spare the details for now and refer to Definition~\ref{def:spiraling} and Theorem~\ref{thm:weak-limits}.



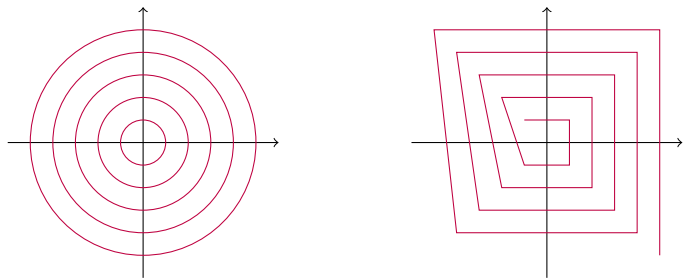
\begin{figure}[ht]
\begin{center}
\begin{tikzpicture}[scale=0.6]
\draw [->] (-3,0)--(3,0);
\draw [->] (0,-3)--(0,3);
\draw [purple] (0,0) circle [radius=0.5];
\draw [purple] (0,0) circle [radius=2];
\draw [purple] (0,0) circle [radius=1];
\draw [purple] (0,0) circle [radius=2.5];
\draw [purple] (0,0) circle [radius=1.5];
\end{tikzpicture}
\hspace{1.5cm}
\begin{tikzpicture}[scale=0.6]
\draw [->] (-3,0)--(3,0);
\draw [->] (0,-3)--(0,3);
\draw [purple] (-0.5,0.5) -- (0.5,0.5)--(0.5,-0.5)-- (-0.5,-0.5)--(-1,1)--(1,1)--(1,-1)--(-1,-1)--(-1.5,1.5)--(1.5,1.5)--(1.5,-1.5)--(-1.5,-1.5)--(-2,2)--(2,2)--(2,2)--(2,-2)--(-2,-2)--(-2.5,2.5)--(2.5,2.5)--(2.5,-2.5);
\end{tikzpicture}
\end{center}
\caption{Two examples of spiraling curves. \emph{Left:} The set of concentric circles $O_\eta$. \emph{Right:} The path $\mathcal{S}_\eta(z_1,...,z_4)$ generated by the points $\{(-a,a),(a,a),(a,-a),(-a,-a)\}$. }
\end{figure}

 Beurling's theory of weak limits \cite{beurl66,beu89,seip04} states an equivalence between sampling sets $\Gamma$ for the Paley-Wiener space $PW_2(\Omega)$ and the uniqueness property for all  weak limits  of translates of $\Gamma$ on the Bernstein space $PW_\infty(\Omega)$. This is a very     powerful result  as the uniqueness property is often easier to verify. The price  one has to pay is that uniqueness has to be shown on a larger space and for all weak limits of translates of the original set. A similar result for sampling measures for the modulation space $M^p(\R)$ was shown by Ascensi \cite{Ascensi}. In case of the measure $\chi_\Gamma \ud \H^1$,  this result states that $\Gamma$ is a Gabor sampling trajectory for $M^p(\R)$ if and only if all weak limits of translates of $\Gamma$ are uniqueness sets for the weighted modulation space   $M^p_{1/\vartheta_s}(\R)$,  for some $s>2$.

In most applications, Beurling's theory is used to derive   necessary conditions on sampling sets. In \cite{janero20} for example, the authors showed that certain weak limits of translates of spiraling curves cannot be uniqueness sets for $PW_\infty(\Omega)$ if the separation exceeds a certain threshold. 
Together with a result from \cite{bewu00}, this yields a full  characterization which  sets of concentric circles and which Archimedes spirals are sampling trajectories in terms of their separation parameter \cite[Theorem A]{janero20}. 

We take a different approach in this paper in that we first  characterize the set of weak limits of translates of spiraling curves and then   show (Theorem~\ref{thm:weak-limits}) that, for windows taken from  the linear span of the Hermite functions, the uniqueness property is automatically satisfied for all   weak limits of $\Gamma$ except for $\Gamma$ itself (and its finite shifts). It then follows that $\Gamma$ is a Gabor sampling trajectory for $M^p(\R)$ if and only if $\Gamma$ is a uniqueness set for $M^p_{1/\vartheta_s}(\R)$  for some $s>2$ (Theorem~\ref{thm:sampling=uniqueness}). Moreover, we prove that certain spiraling curves are indeed uniqueness sets.
In particular, our main result is the following.

\begin{theorem}\label{thm:main}
Let $g\in \text{span}\{h_n:   n\in\N_0\}$, $1\leq p<\infty$,  and $\eta>0$. Moreover, let \begin{enumerate}
\item[(a)] $P $ be  a star shaped polygon in $\R^2$ such that $0\in \ker (P)$ and no line passing through an edge of the polygon meets the origin,
\item[(b)] $ z_1,...,z_n \in\R^2,\ n\geq 3$,  be such that $ \text{arg}(z_1),...,\text{arg}(z_n),\text{arg}(z_1 )$ is strictly increasing or strictly decreasing modulo $2\pi$. 
\end{enumerate}
The following sets are Gabor sampling trajectories for $M^p(\R)$:
\begin{itemize}\item[(i)] the set of concentric circles
$
\mathcal{O}_\eta:=\{(x,y)\in\R^2 :\ x^2+y^2=\eta^2k^2,\ k\in \N\},
$
\item[(ii)] the collection of star shaped polygons
$
\mathcal{P}_\eta:=\{(x,y)\in\R^2:\ (x,y)\in \eta k P,\ k\in \N\},
$
\item[(iii)] $\mathcal{S}_\eta(z_1,...,z_n)$, the path generated by the sequence of vectors
$$
\{\eta z_1,\eta z_2,...,\eta z_n,2\eta z_1,2\eta z_2,...,2\eta z_n,3\eta z_1,3\eta z_2,...\}.
$$
\end{itemize}
\end{theorem}


\noindent This paper is organized as follows. After a section dedicated to preliminaries, we study some basic examples of windows and trajectory sets in Section~\ref{sec:nec} and give necessary and sufficient conditions for Gabor sampling trajectories in that setting. Then in Section~\ref{sec:spiraling} we introduce spiraling curves and study their weak limits of translates which then allows us to prove our main result.

\section{Preliminaries}
\subsection{Notation}

Throughout this paper we adopt the following conventions. 
We write $\vec{d}$ for vectors in $\mathbb{S}^1$ and 
$\vec{\ell}(\theta):=(\cos(2\pi \theta),\sin(2\pi\theta))$, $\theta\in\mathbb{T}\cong [0,1)$ as well as $\vec{e}_i$ for the standard basis vectors in $\R^2$. For $\vec{d}\in \mathbb{S}^1$ we write $\vec{d}_\bot$ for the vector obtained by rotating $\vec{d}$ by $\pi/2$. The space of continuous, compactly supported functions is denoted by $C_c(\R^d)$, and 
we write 
$A\lesssim B$  if there exist $C>0$ such that $A\leq CB$.

\subsection{Hausdorff measure and trajectory sets}
Let $E\subset\R^2$. The \emph{$1$-dimensional Hausdorff measure} of $E$ is given by 
$$
\H^1(E)=\lim_{\delta\searrow 0}\inf\left\{2\sum_j r_j:\ E\subset \bigcup_j B(x_j,r_j)\ \text{and}\ r_j\leq \delta\right\}.
$$
When restricted to a line,   $ \H^1$ equals the $1$-dimensional Lebesgue measure.

Let $\varphi:[0,1)\rightarrow[0,\infty)$ be continuous in $0$. A (locally finite Borel) measure $\mu$ is called \emph{$\varphi$-regular} if for every $z\in\R^2$ and every $R\in (0,1)$  one has
$$
\mu(B_R(z))\leq \pi\varphi(R) R.
$$
 A set $E\subset \R^2$ is called \emph{$\varphi$-regular} if $\H^{1}|_E$ is $\varphi$-regular.

 A \emph{trajectory} $\Gamma\subset \R^d$ is the image of a curve $\gamma:\R\rightarrow\R^d$, i.e. $\Gamma=\gamma(\R)$, such that each restriction of $\gamma$ to a finite interval is rectifiable. A \emph{trajectory set} $\Gamma$ is   a countable collection of trajectories.

\subsection{The short-time Fourier transform and modulation spaces}
Let $z=(x,\xi)\in \R^2$. A \emph{time-frequency shift}  of a function $g$ is given by
$$
\pi(z)g(t)=M_\xi T_x g(t)=e^{2\pi i\xi t}g(t-x).
$$
Given a window $g\in L^2(\R)$, the \emph{short-time Fourier transform} of $f\in L^2(\R)$ is defined as
$$
V_g f(z)=\langle f,\pi(z)g\rangle =\int_\R f(t)\overline{g(t-x)}e^{-2\pi i \xi t}\ud t.
$$
The short-time Fourier transform satisfies the following \emph{orthogonality relation}
$$
\int_{\R^2}V_{g_1}f_1(z)\overline{V_{g_2}f_2(z)}\ud z=\langle f_1,f_2\rangle \overline{\langle g_1,g_2\rangle},
$$
which implies that, if $\|g\|=1$, $V_g:L^2(\R)\rightarrow L^2(\R^2)$ is an isometry.
 
Let $\vartheta:\R^2\rightarrow \R^+$ be a \emph{submultiplicative weight function}, i.e. $\vartheta(z+w)\leq \vartheta(z)\vartheta(w)$. The 
\emph{(weighted) modulation spaces}  $M_{\vartheta }^p(\R)$, $1\leq p\leq\infty$, can be defined by
$$
M_{\vartheta }^p(\R):=\left\{f\in\mathcal{S}'(\R):\ V_{h_0}f \cdot\vartheta \in L^p(\R)\right\},
$$
where $h_0$ denotes  the normalized Gaussian and $\mathcal{S}'(\R)$ the space of tempered distributions. Modulation spaces are Banach spaces when equipped with the natural norm $\|f\|_{M_\vartheta} =\|V_{h_0}f\cdot\vartheta\|_{L^p}$. See for example \cite{groe1}   for a thorough introduction to modulation spaces. Throughout this paper we   only consider  \emph{polynomial weight functions}, i.e.,
$$
\vartheta_s(z)=(1+|z|)^s,\quad z\in\R^2,\ s\geq 0.
$$
We note here that $M_{\vartheta_s}^1(\R)$ is closed under pointwise multiplication with functions from the \emph{weighted Fourier algebra} $A_{\vartheta_s}(\R)$ 
$$
A_{\vartheta_s}(\R):=\left\{f\in C_0(\R):\ \int_\R|\widehat{f}(\xi)|\vartheta_s(\xi)\ud\xi<\infty\right\}.
$$
In particular, for $f\in A_{\vartheta_s}(\R)$ and $g\in M_{\vartheta_s}^1(\R)$ we have
\begin{equation}\label{eq:fourier-algebra}
\|fg\|_{M_{\vartheta_s}^1 }\leq \|f\|_{A_{\vartheta_s}}\|g\|_{M_{\vartheta_s}^1 },
\end{equation}
see for example the arguments in \cite[Proposition 4.13]{ja18} which can easily be adapted to the weighted case.

Let $g\in M^1(\R)$. We call a trajectory set $\Gamma\subset\R^2$   a \emph{Gabor sampling trajectory} for $M^p(\R)$ if there exist constants $A,B>0$ such that 
\begin{equation}\label{eq:def-gabor-trajectory}
A\|f\|^p_{M^p}\leq \int_\Gamma |V_g f(z)|^p\ud \H^1(z)\leq B\|f\|^p_{M^p},\quad f\in M^p(\R).
\end{equation}
If only the upper bound is satisfied, then we call $\Gamma$ a \emph{Gabor Bessel trajectory}.
If $p=2$, then $M^p(\R)=L^2(\R)$, and \eqref{eq:def-gabor-trajectory} is equivalent to $\{\pi(z)g\}_{z\in\Gamma}$ forming a continuous frame, see \cite{alanga93,ranade06}.

The trajectory set $\Gamma$ is called a \emph{uniqueness set} for $M^p(\R)$, if  $V_gf|_\Gamma=0$, $f\in M^p(\R)$,  implies $f=0$.

\subsection{Hermite windows and polyanalytic functions}

A function $F:\C\rightarrow \C$ is called \emph{polyanalytic of order $n$} if it satisfies the higher order Cauchy-Riemann equation $(\bar\partial)^{n+1}F=0$. In that case, $F$ can be written as 
\begin{equation}\label{poly1}
F(z)=F(z,\overline{z})=\sum_{k=0}^n F_k(z)\overline{z}^k,\quad z\in\C
\end{equation}
where $F_0,\ldots,F_n:\C\rightarrow\C$ are holomorphic functions.

The \emph{Hermite functions} are given by
\begin{equation}\label{eq:def-hermite}
h_n(t)=\frac{2^{1/4}}{\sqrt{n!}}\left(\frac{-1}{2\sqrt{\pi}}\right)^n e^{\pi t^2}\frac{d^n}{dt^n}\left(e^{-2\pi t^2}\right),\quad n\in\N_0.
\end{equation}
For $f\in M^p_{\vartheta_s}(\R)$, $1\leq p<\infty$, the function $F$ given by
\begin{equation}\label{def-poly-bargmann-trafo}
F(z)=V_{h_n}f(\overline{z})e^{\pi (z^2-\overline{z}^2)/4} e^{\pi |z|^2/2},
\end{equation}
is   polyanalytic   of order $n$ \cite{abgroe12},
where we identify  $z=(x,\xi)\in\R^2$ with $z=x+i\xi\in\C$. In particular, if $g=\sum_{k=0}^n \alpha_k h_k$, then $V_{g}f(\overline{z})e^{\pi (z^2-\overline{z}^2)/4} e^{\pi |z|^2/2}$ is polyanalytic of order $n$.

A polyanalytic function $F$ of order $n$ is called \emph{reduced} if it can be written as
$$
F(z)=\sum_{k=0}^n F_k(z)|z|^{2k},\quad F_k\text{ holomorphic.}
$$
Reduced polyanalytic functions  satisfy the following  Cauchy-type formula \cite[Section 1.3, (11)]{balk97}. For a similar Cauchy-type formula for true polyanalytic functions, see \cite{abspe18-sieve}.
\begin{lemma}[Balk]
Let $F$ be a reduced polyanalytic function of order $n$ in $B_R(0)$, $0<R_0<R_1<\ldots<R_n<R$, and let $\Gamma_k:=\{z:\ |z|=R_k\}$. For every $z\in B_{R_0}(0)$
\begin{equation}\label{eq:cauchy-formula}
F(z)=\frac{1}{2\pi i}\sum_{k=0}^n P_k(|z|^2)\int_{\Gamma_k}\frac{F(t)}{t-z}dt,
\end{equation}
where  
$
P_k(t):=\prod_{j\neq k}\frac{R_j^2-t}{R_j^2-R_k^2}.
$
\end{lemma} 


\subsection{The metaplectic rotation}
 Let us  denote the  rotation matrices in $\R^2$ by 
$
R(\theta)=\left(\begin{smallmatrix}
\cos( 2\pi\theta) & -\sin (2\pi\theta)\\ \sin (2\pi\theta) &\cos (2\pi\theta)
\end{smallmatrix}\right),\ \theta\in\mathbb{T}.
$
The \emph{metaplectic rotation} of $f\in L^2(\R)$ is given in terms of the Hermite basis $\{h_n\}_{n\in\N_0}$
$$
\mu(\theta) f:=\sum_{n\geq 0}e^{-2\pi in\theta}\langle f,h_n\rangle h_n.
$$
Clearly, $\mu(\theta)$ is a unitary operator on $L^2(\R)$ with $\mu(\theta)^\ast=\mu(-\theta)$, and  $\mu(\theta) h_n=e^{-2\pi in\theta} h_n$.
For $f,g\in L^2(\R)$, the  standard rotation of the argument of the short-time Fourier transform and the metaplectic rotation are connected via the formula
\begin{equation}\label{eq:metaplectic}
V_gf (R(\theta) z)=e^{\pi i (x\omega-x'\omega')}V_{\mu({ \theta}) g}\mu({ \theta})f(z),\quad z=(x,\omega),\text{ and }(x'\omega')=R(\theta)  z,
\end{equation}
or equivalently via $\pi(R(\theta)  z)=e^{-\pi i (x\omega-x'\omega')}\mu(-\theta)\pi(z)\mu(\theta)$.
This is a special case of the symplectic covariance of the Schr\"odinger representation, see \cite[Chapters 1 \& 2]{fo89}, and \cite[Chapter 9]{groe1}. 

As $\vartheta_s$ is radially symmetric, \eqref{eq:metaplectic} implies that $\|f\|_{M^1_{\vartheta_s} }=\|\mu(\theta)f\|_{M^1_{\vartheta_s} }$, which allows   to extend the operator $\mu(\theta)$ to $M^\infty_{1/\vartheta_s}(\R)=M^1_{\vartheta_s}(\R)^\ast$ via
$$
\langle \mu(\theta)f,g\rangle _{M^\infty_{1/\vartheta_s}\times M^1_{\vartheta_s}}=\langle  f,\mu(-\theta)g\rangle _{M^\infty_{1/\vartheta_s}\times M^1_{\vartheta_s}},\quad f\in M^\infty_{1/\vartheta_s}(\R),\ g\in M^1_{\vartheta_s}(\R).
$$
Consequently,  \eqref{eq:metaplectic} remains valid for $f\in M^\infty_{1/\vartheta_s}(\R),$ and $g\in M^1_{\vartheta_s}(\R)$.

 \section{General trajectory sets}\label{sec:nec}

  In this section, we present some basic necessary and sufficient conditions for  Gabor sampling trajectories. Moreover, we study the particular cases  of (i)  Gaussian window with general trajectory sets and (ii)   sampling on parallel lines for general windows.

\subsection{Relative density}\label{subsec:nec}

 We call a trajectory set $\Gamma\subset\R^2$ \emph{$(m,R)$-dense} if
 $$
 \inf_{z\in\R^2}\H^1(\Gamma\cap B_R(z) )\geq m>0, 
 $$
 and \emph{relatively dense} if there exist constants $m,R>0$ such that $\Gamma$  is $(m,R)$-dense.
It turns out that relative density is  a necessary condition for $\Gamma$ being a Gabor sampling trajectory for $M^p$, see \cite[Theorem 8 \& 10]{Ascensi}:

\begin{proposition}[Ascensi]
Let $g\in M^1(\R)$ and $1\leq p<\infty$.
If $\Gamma\subset \R^2$ is a Gabor Bessel trajectory for $M^p(\R)$, then there exist $M,R>0$ such that 
\begin{equation}\label{eq:nec-bessel}
\sup_{z\in\R^2}\H^1(\Gamma\cap B_{R}(z))\leq M.
\end{equation}
 If $\Gamma$ is a Gabor sampling trajecory for $M^p(\R)$, then $\Gamma$ is relatively dense.
\end{proposition}

\noindent Note that the proof of \eqref{eq:nec-bessel} can easily be adapted to the case $p=2$ and general $g\in L^2(\R)$ and   explicit   upper bounds of the Bessel constant can be derived for example from \cite{abspe17-sampta,abspe18-sieve}. 

The question whether relative density  is also necessary for general windows and Gabor sampling trajectories for $L^2(\R)$ remains open. Drawing comparison to the discrete \cite{rast95} and planar cases \cite{jasp20} however suggests that relative density should indeed be a necessary requirement.

\subsection{Gaussian window}\label{subsec:gauss}

In \cite{orce98}, Ortega-Cerd{\`a}   fully characterized the  sampling measures for the Bargmann-Fock space of entire functions which corresponds to the short-time Fourier transform with Gaussian window. The result goes as follows.

\begin{theorem}[Ortega-Cerd{\`a}]\label{thm:ortega-cerda}
Let $1\leq p<\infty$. The measure $\mu$ is a sampling measure for $V_{h_0}(M^p(\R))$ if and only if there exist constants $R,\delta,M>0$ and $N\in\N$ such that 
$$
(i)\  \sup_{z\in\R^2}\mu(B_R(z))\leq M,\quad \quad\quad
 (ii)\ \inf_{z\in\R^2}\frac{n(R,N,\delta,z)}{R^2}>1,
 $$
 where 
$n(R,N,\delta,z)$ is calculated by the following rule: take $S_R(z):=z+[-R/2,R/2)^2$   and cover it with $N^2$ smaller squares of sidelength $R/N$. Then $n(R,N,\delta,z)$ denotes the number of   smaller squares $s\subset S_R(z)$ for which $\mu(s)\geq \delta$.
\end{theorem}

\noindent When specifying the measure $\mu$ to be the $1$-dimensional Hausdorff measure on $\Gamma$, we subsequently show that, with a minor extra condition on the trajectory set $\Gamma$,   condition (ii) is equivalent to $\Gamma$ being relatively dense.

\begin{corollary}
Let $h_0$ be the standard Gaussian window and $\Gamma\subset \R^2$ be a $\varphi$-regular trajectory set. Then $\Gamma$ is a Gabor sampling trajectory   for $M^p(\R), \ 1\leq p<\infty,$ if and only if there exist $R,m,M>0$ such that 
\begin{equation}\label{eq:cor-gauss}
m\leq \H^1(\Gamma\cap B_R(z))\leq M,\quad \text{for every }z\in\R^2.
\end{equation}
\end{corollary}
\proof The second inequality of \eqref{eq:cor-gauss} is condition (i) in Theorem \ref{thm:ortega-cerda}. Hence, we have to show the equivalence of the left hand side inequality and (ii). If (ii) holds, then
$$
\H^1(\Gamma\cap B_{R}(z))\geq \H^1(\Gamma\cap S_R(z))\geq \delta n(R,N,\delta,z)>\delta R^2.
$$
Now let $\Gamma$ be relatively dense and choose $N$ large enough such that $\varphi(\sqrt{2}R/N)\leq \sqrt{2}\varphi(0)$. Then, using the lower bound in \eqref{eq:cor-gauss} and the $\varphi$-regularity of $\Gamma$, we get
\begin{align*}
m&\leq \H^1(\Gamma\cap B_R(z)) \leq \H^1(\Gamma\cap S_{2R}(z))
\\ 
&\leq   n(2R,N,\delta,z)\sup_{w\in S_{2R}(z)}\H^1(\Gamma\cap S_{\frac{2R}{N}}(w))  +(N^2-n(2R,N,\delta,z))\delta
\\
&\leq  n(2R,N,\delta,z)2\pi\varphi(0)\frac{ R}{N}+(N^2-n(2R,N,\delta,z))\delta,
\end{align*}
for some constant $C>0$.
Regrouping the terms yields
$$
n(2R,N,\delta,z))\geq N\frac{m-N^2\delta}{2\pi R-\delta N}.
$$
Choosing $\delta=\frac{1}{N^2}\min\{m/2,\pi R\}$ shows that $n(2R,N,\delta,z))$ can be   arbitrarily large as $N\rightarrow\infty$. In particular, there exists, $N\in\N$ such that $(ii)$ is satisfied. \hfill $\Box$

\begin{remark}
In \cite{pr87} it is shown that
 $\varphi$-regularity with $\varphi(0) =1$ can be
considered as a quantitative strengthening of rectifiability. On the
other hand, it is possible to construct fractal sets that are
$\varphi$-regular if $\varphi(0) > 1$, which shows that $\varphi$-regularity for some $\varphi$ is a rather mild assumption. 
\end{remark}

\subsection{The case of    parallel lines}
We now study  sampling and uniqueness properties of sets of parallel lines which follow from simple arguments. The results however will be  useful later when we  study  spiraling curves.
\begin{proposition}\label{prop:parallel-lines-frame}
Let $g\in L^2(\R)$, $\Lambda\subset \R$ be countable, and $\vec{d}=R(\theta)\vec{e}_2$. The collection of parallel lines $L_{\vec{d},\Lambda}$, where 
$$
L_{\vec{d},\Lambda}:= \big\{t\vec{d} +  \lambda\vec{d}_\bot:\ t\in\R,\ \lambda\in\Lambda\big\},
$$
 is a Gabor sampling trajectory for $L^2(\R)$ with sampling bounds $A,B>0$ if and only if
\begin{equation}\label{eq:cond-frames-lines}
A\leq \sum_{\lambda\in\Lambda}|\mu(\theta) {g}(t-\lambda  )|^2\leq B,\quad \text{ for a.e. }t\in\R.
\end{equation} 
\end{proposition}
\proof
First, by \eqref{eq:metaplectic} we can rotate the problem and  assume   $\vec{d}=\vec{e_2}$.
Writing the 
short-time Fourier transform as $V_gf(x,\xi)=\mathcal{F}\big(f \overline{T_xg}\big)(\xi)$ and using Parseval's identity yields
\begin{align*}
\int_{L_{\vec{e_2},\Lambda}}|V_gf(z)|^2\ud\H^1(z)&=\sum_{\lambda\in\Lambda}\int_\R |V_gf( \lambda  ,\xi)|^2\ud\xi=\sum_{\lambda\in\Lambda}\int_\R |f(t)\overline{g(t-\lambda )}|^2\ud t
\\
&=\int_\R |f(t)|^2\left(\sum_{\lambda\in\Lambda}|g(t-\lambda  )|^2\right)\ud t,
\end{align*}
where changing the order of integration and summation is allowed by either \eqref{eq:cond-frames-lines} or the existence of the upper sampling bound .\hfill $\Box$\medskip

\noindent The uniqueness property of parallel lines  on the distribution space $M^\infty_{1/\vartheta_s}(\R)$ can be characterized in a similar fashion.

\begin{proposition}\label{prop:lines-uniqueness}
Let $g\in M^1_{\vartheta_s}(\R)\cap A_{\vartheta_s}(\R)$, $s\geq 0$, $\Lambda\subset \R$ be countable and $\vec{d}=R( \theta)\vec{e}_2$. The collection of parallel lines $L_{\vec{d},\Lambda}$ is a uniqueness set for $M_{1/\vartheta_s}^\infty(\R)$ if and only if 
$$
\bigcup_{\lambda\in\Lambda}\widetilde{\emph{supp}}(T_{\lambda  }\mu(\theta)g)=\R,
$$
where $\widetilde{\emph{supp}}(g):=\{t\in\R:\ g(t)\neq 0\}$ denotes the \emph{effective support} of $g$.
\end{proposition}

\begin{remark}
In \cite[Lemma 27]{Ascensi}, a similar result was shown for planar uniqueness sets. Since some technical details are left out there, we   decided to  include the proof here.
\end{remark}

\noindent \proof For $f\in M^\infty_{1/\vartheta_s}(\R)$, $g\in M^1_{\vartheta_s}(\R)\cap A_{\vartheta_s}(\R)$ and $h\in M^1_{\vartheta_s}(\R) $,  we define $f g$ via   $\langle f g,h\rangle_{M^\infty_{1/\vartheta_s}\times M^1_{\vartheta_s}}$ $=\langle f,gh\rangle_{M^\infty_{1/\vartheta_s}\times M^1_{\vartheta_s}}$. By \eqref{eq:fourier-algebra}, it then follows that $f g\in M^\infty_{1/\vartheta_s}(\R)$ as
$$
|\langle f g,h\rangle_{M^\infty_{1/\vartheta_s}\times M^1_{\vartheta_s}}|=|\langle f,gh\rangle_{M^\infty_{1/\vartheta_s}\times M^1_{\vartheta_s}}|\leq \|f\|_{M^\infty_{1/\vartheta_s}}\|gh\|_{M^1_{\vartheta_s}}\leq \|f\|_{M^\infty_{1/\vartheta_s}}\|g\|_{A_{\vartheta_s}}\| h\|_{M^1_{\vartheta_s}}.
$$ 
Moreover, since $A_{\vartheta_s}(\R)$ is invariant under translations it follows that
$f\overline{T_xg}\in M^\infty_{1/\vartheta_s}\subset \mathcal{S}'(\R)$ and consequently, $\mathcal{F}\big(f\overline{T_xg}\big)\in\mathcal{S}'(\R)$ is   a well defined tempered distribution.

Since $M^1_{\vartheta_s}(\R)$ is weak-$\ast$ dense in $M^\infty_{1/\vartheta_s}(\R)$, and since the short-time Fourier transform of a distribution in $M^\infty_{1/\vartheta_s}(\R)$ is   continuous   in $\R^2$, it follows that $V_g f (x,\xi)=\mathcal{F}\big(f\overline{ T_xg}\big)(\xi)$.

As before we may rotate the problem and assume  $\vec{d}=\vec{e_2}$. Now,   $V_gf|_{L_{\vec{e_2},\Lambda}}=0$ if and only if the distributions $f\overline{T_{\lambda  }g},\ \lambda\in\Lambda$, are zero. This in turn is equivalent to the support of $f$ and the effective support of $T_{\lambda  }g $ being  disjoint for every $\lambda\in\Lambda$. 
\hfill $\Box$

\subsection{Connection to discrete sampling}

 Proposition~\ref{prop:parallel-lines-frame} shows that $\Gamma$ being relatively dense is not sufficient for $\Gamma$ to be a Gabor sampling trajectory. A natural question is therefore whether for every  $g\in L^2(\R)\backslash\{0\}$ there exists $R^\ast>0$  such that every $(\gamma,R)$-dense trajectory set $\Gamma$ is a Gabor sampling trajectory  if $R\leq R^\ast$.
We   follow the approach of \cite{jasp20} to show that such  $R^\ast$ does in fact exist for a certain class of window functions.


Let us write $Q_R(z):=Rz+ [-R/2,R/2)^2$ and recall the definition of the Sobolev space  
$
 H^1(\R):=\Big\{f\in L^2(\R):\ \int_\R (1+|\xi|^2)|\widehat{f}(\xi)|^2\ud \xi<\infty\Big\}.
$

\begin{proposition}\label{thm:planar-sampling-frames}
Let $g,tg\in H^1(\R)$, $R>0$ be chosen such that 
 \begin{equation}\label{eq:con-frame-suzhou}
 \Delta:=\frac{2R}{\pi}\left(\|g'\|_2+\|tg\|_2+\frac{2R}{\pi}\|tg'\|_2\right)<\|g\|_2,
 \end{equation} and $\Gamma\subset \R^2$ be a trajectory set. 
  If  
$$
0<m\leq \H^1(\Gamma\cap Q_R(z))\leq M<\infty,\quad \text{for every }z\in\R^2,
$$ 
then 
$$
m\big(\|g\|_2-\Delta\big)^2\|f\|^2\leq \int_\Gamma |V_gf (z)|^2 \ud\H^1(z)\leq M\big(\|g\|_2+\Delta\big)^2\|f\|^2.
$$
\end{proposition}
\proof It is shown in \cite{suzho02} that for the particular choice of $g$ and  $R$, taking arbitrary points $z_{n}\in Q_R(n),\ n \in\Z^2$,   yields a discrete frame $\{\pi(z_{n})g\}_{n\in\Z^2}$ for $L^2(\R)$ with   uniform frame bounds $A=\big(\|g\|_2-\Delta\big)^2$ and $B=\big(\|g\|_2+\Delta\big)^2$. 

 For every $n \in\Z^2$ there exists $z_{n }\in \Gamma\cap Q_R({n})$ such that
$$
|V_g f(z_{n})|^2\geq \frac{1}{\H^1(\Gamma\cap Q_R({n}))}\int_{\Gamma\cap Q_R({n})}|V_g f(z)|^2\ud\H^1(z).
$$
Then, as every choice of points $z_{n }\in Q_R({n})$ generates a Gabor frame with uniform upper bound $B$,  we have
\begin{align*}
\int_\Gamma |V_{g}f(z)|^2\ud\H^1(z)&=\sum_{n\in\Z^2}\int_{\Gamma\cap Q_R({n})} |V_{g}f(z)|^2\ud\H^1(z)
\\
&\leq \sum_{n\in\Z^2}\H^1(\Gamma\cap Q_R({n})) |V_{g}f(z_{n})|^2\leq MB \|f\|^2.
\end{align*}
The lower sampling bound follows with a similar argument. \hfill $\Box$

\begin{remark}  
This construction works for general   measures and gives a characterization of sampling measures for this class of window functions. 

\end{remark}

\section{Spiraling curves  }\label{sec:spiraling}
\subsection{Weak limits}

There are multiple ways of defining weak limits of trajectory sets. The definition in \cite{janero20} for example adapts the initial definition by Beurling \cite{beu89} which is given in terms of a geometric condition. Here, it will be convenient to work with a stronger notion that was introduced to define weak limits of measures, see \cite{Ascensi}.
\begin{definition}
Let $\Gamma,\Gamma',\Gamma_n$ be  trajectory sets in $\R^2$.
We say that $\{\Gamma_n\}_{n\in\N}$ \emph{converges weakly} to $\Gamma$ if  
$$
\int_{\Gamma_n}\phi(z) \ud\H^1(z)\rightarrow\int_{\Gamma}\phi(z) \ud\H^1(z),
$$
for every nonnegative function $\phi\in C_c(\R^2)$. In that case we write $\Gamma_n\stackrel{w}{\rightarrow}\Gamma$.

We say that $\Gamma'$ is a weak limit of translates of $\Gamma$ if there exists a sequence $\{z_n\}_{n\in\N}$ such that $z_n+\Gamma\stackrel{w}{\rightarrow}\Gamma'$ and define $\mathcal{W}_\Gamma$ as the set of all weak limits of translates of $\Gamma$.
\end{definition}

\noindent The following characterization of Gabor sampling trajectories is an immediate consequence of a characterization of samping measures for the short-time Fourier transform given by Ascensi \cite[Theorem 14]{Ascensi}. One can think of this result as a time-frequency analog of the classical result by Beurling \cite[Theorem 3, pg. 345]{beu89} where the Gabor sampling property on $M^p(\R)$ is connected to the uniqueness   property  of all weak limits of translates on the larger space $M^p_{1/\vartheta_s}(\R)$. Note that Ascensi's result requires a special class of  windows   which we define in a simplified version taylored to polynomial weights.
\begin{definition}
We say that  $1/\vartheta_s$ \emph{controls} $V_gg$ if
there exist a  function $d:\R^+\rightarrow\R^+$ which is decreasing with $d(r)\rightarrow 0$, as $r\rightarrow \infty$, and satisfies
$
|V_gg(z)|\leq  d(|z|)/\vartheta_s(z).
$
The class $\mathcal{M}^1_s(\R)$, $s>2$, is then given by 
$$
\mathcal{M}^1_s(\R):=\big\{g\in L^2(\R):\  1/\vartheta_s  \text{ controls }V_gg\big\}.
$$
\end{definition}

\begin{theorem}[Ascensi]\label{thm:ascensi}
Let $g\in\mathcal{M}^1_s(\R)$, for some $s>2$, and $1\leq p<\infty$. Then $\Gamma$ is a Gabor sampling trajectory  for $M^p(\R)$ if and only if every $\Gamma'\in \mathcal{W}_\Gamma$ is a set of uniqueness for $M^p_{1/\vartheta_s}(\R)$.
\end{theorem}

\subsection{Spiraling curves and their weak limits}

The notion of a spiraling curve was introduced in \cite{janero20}. This class of trajectory sets  includes a wide range of natural examples such as the concentric circles or the Archimedes spiral. 
In this paper, we use a slightly more restrictive notion of spiraling curves  that still includes the main examples from  \cite{janero20} and that allows a full characterization of the set of weak limits of translates.

\begin{definition}[Spiraling curve]\label{def:spiraling} Let $\mathcal{I}_\Gamma\subset \mathbb{T} $ be a finite set.
A   trajectory  set $\Gamma$ is called \emph{spiraling} if the following list of  conditions is satisfied:

\noindent   If $\beta\in \mathbb{T}\backslash \mathcal{I}_\Gamma$, then
\begin{enumerate}
\item[(A.i)]  \emph{(Escape Cone)} for $\alpha:=\min\{\emph{dist}(\beta,\mathcal{I}_\Gamma),1/8\}$   the intersection of $\Gamma$ with the cone
$$
\mathcal{C}_{\alpha,\beta}:=\{(r \cos 2\pi \theta,r \sin 2\pi\theta):\ r\geq0,\ \beta-\alpha\leq \theta\leq \beta+\alpha\}
$$
can be parametrized in polar coordinates as
$$
\gamma_\beta(\theta)=(r_\beta(\theta)\cos 2\pi \theta,r_\beta(\theta)\sin 2\pi \theta)
$$
with $\theta\in\bigcup_{k\in\N} [k+\beta-\alpha,k+\beta+\alpha]$ and $r_\beta$ a nonnegative $C^2$-function on each interval whose $C^2$-norm is globally bounded.


\item[(A.ii)]  \emph{(Asymptotic radial monotonicity)} there exist $K\in\N$ such that for any $\theta\in[ -\alpha, \alpha] $ fixed, the sequence $r_\beta(k+\beta+\theta)$ is strictly increasing for $k\geq K$.

\item[(A.iii)] \emph{(Asymptotic flatness)} the curvature of $\gamma_\beta$, denoted by $\kappa_\beta$, tends to $0$ as its input goes to infinity. To be more precise, we assume 
\begin{equation}\label{eq:as-flat}
\sup_{\theta\in (-\alpha,\alpha)} \kappa_\beta(k+\beta+\theta)\rightarrow 0,\quad k\rightarrow \infty.
\end{equation}


\item[(A.iv)] \emph{(Asymptotic equispacing)} there exist continuous functions $\eta_\beta:[-\alpha,\alpha]\rightarrow \R^+,\rho_\beta:[-\alpha,\alpha]\rightarrow \R$ such that   
$$
\sup_{\theta\in[-\alpha,\alpha]}\left|r_{ \beta}(k+\beta+\theta)-\eta_\beta(\theta)k -\rho_\beta(\theta)\right|\rightarrow 0,\quad k\rightarrow\infty.
$$
\item[(A.v)] \emph{(Asymptotic velocity)} there exists a continuous function $\vec{d}_\beta:[ -\alpha, \alpha]\rightarrow\mathbb{S}^1$ such that $\vec{d}_\beta(\theta)$ is  non-collinear with $\vec{\ell}(\theta)$  and
\begin{equation}\label{eq:as-velo}
\sup_{\theta\in[-\alpha,\alpha]}\left\|\frac{\gamma_\beta^\prime(k+\beta+\theta)}{\|\gamma_\beta^\prime(k+\beta+\theta)\|}-\vec{d}_\beta(\theta)\right\|\rightarrow 0,\quad k\rightarrow\infty.
\end{equation}
\end{enumerate}

\noindent 
If $\beta\in \mathcal{I}_\Gamma$, then $(A.i)-(A.v)$ hold with the following modifications
\begin{enumerate}
\item[(B.i)] for $\alpha:=\min\{\emph{dist}(\beta,\mathcal{I}_\Gamma\backslash\{\beta\}),1/8\}$, the parametrization $\gamma_\beta|_{[k+\beta-\alpha,k+\beta+\alpha]}$ is  continuous and the restrictions of $\gamma_\beta$ to   ${(k+\beta,k+\beta+\alpha]}$ and ${[k+\beta-\alpha,k+\beta)}$ are $C^2$-functions with bounded $C^2$-norms.  
\item[(B.iii)] The supremum in \eqref{eq:as-flat} is replaced by an supremum over $(-\alpha,\alpha)\backslash\{0\}$.
\item[(B.v)]  The velocity vector $ \vec{d}_\beta  $ is continuous on $[ -\alpha,\alpha]\backslash\{0\}$. Morover,   $\lim_{\theta \nearrow 0}\vec{d}_\beta(\theta)$ and $\lim_{\theta {\searrow }0}\vec{d}_\beta(\theta)$ exist, and the supremum in \eqref{eq:as-velo} is taken over $[-\alpha,\alpha]\backslash\{0\}$.
\end{enumerate}
\end{definition}


\begin{remark} 


(i) The original definition of spiraling curves  \cite[Section 3.3]{janero20} only assumed that there exist at least one angle such that the conditions (A.i)-(A.v) are satisfied.
Let us mention here, that our additional assumptions are also met by the examples mentioned in \cite{janero20}. In particular, the set of concentric circles is also a spiraling curve in the sense of this paper.

(ii) Since the collection of star shaped polygons and paths generated by a set of points consist of line segments (mostly parallel and equispaces within escape cones), it is a straightforward task to show that these trajectory sets are indeed spiraling curves. The only escape cone that needs more attention is the one that contains the line segments $s({kz_n, (k+1)z_1})$. In the limit however, these are parallel and equispaced line segments.
 
\end{remark}

\noindent Subsequently, we give a full characterization of the set of all weak limits of translates of spiraling curves.
To do so, we establish two technical lemmas that describe the weak limits of $z_k+\Gamma$ according to   a certain geometric condition  on the sequence $\{z_k\}_{k\in\N}$.

\begin{lemma}\label{lem:lines1}
Let $z_k =  - r_k \vec{\ell}(\theta_k),\ k\in\N,$ be  an unbounded sequence. If there exist no pair $(\mathcal{N},\gamma)$, $\mathcal{N}\subset \N,\ \gamma\in \mathbb{T}$,  such that $$(i)\ \{z_k\}_{k\in\mathcal{N}}\text{ is unbounded \quad and \quad } (ii)\ \{r_{k}\sin(2\pi(\theta_{k}-\gamma))\}_{k \in\mathcal{N}}\text{ is bounded,}$$ then there exist a subsequence $\{z_{k_n}\}_{n\in\N}$, an angle $\beta\in\mathbb{T}\backslash\mathcal{I}_\Gamma$, and a constant $\tau\in\R$,  such that  (a) 
 $-z_{k_n}\in   (\mathcal{C}_{\alpha,\beta})^\circ$, for every $n\in\N$, (b)
 $\theta_{k_n}\rightarrow \theta^\ast$, $n\rightarrow\infty$,
and
$$
(c)\ z_{k_n}+\Gamma\stackrel{w}{\rightarrow}\tau\vec{d}_\bot+L_{\vec{d},\lambda\Z},
$$
where
$\vec{d}=\lim_{n\rightarrow\infty} \vec{d}_\beta(\theta_{k_n}-\beta)$ and $\lambda= \eta_\beta(\theta^\ast-\beta) \sin\big(\arccos(\vec{\ell}(\theta^\ast) \cdot \vec{d}\ )\big).$

Moreover, for every $\tau\in\R$, $\beta\in \mathbb{T}\backslash \mathcal{I}_\Gamma$ and for  $\vec{d}=\lim_{\theta\searrow 0}\vec{d}_{\beta}(\theta)$ as well as for  $\vec{d}=\lim_{\theta\nearrow 0}\vec{d}_{\beta}(\theta)$, there exist   a sequence $\{z_k\}_{k\in\N}$ such that 
$$z_k+\Gamma\stackrel{w}{\rightarrow}\tau\vec{d}_\bot+L_{\vec{d},\lambda\Z},$$ with  $\lambda= \eta_\beta(0) \sin\big(\arccos(\vec{\ell}(\beta)\cdot \vec{d})\big)$.
\end{lemma}
\begin{remark}
 The distinction of the two directions (which coincide for $\beta\in\mathbb{T}\backslash\mathcal{I}_\Gamma$) in   Lemma~\ref{lem:lines1} is necessary as, for $\beta\in\mathcal{I}_\Gamma$, the sets of parallel lines  with the left and right limits of the asymptotic velocity as directions are both included in  $\mathcal{W}_\Gamma$.

\end{remark}

\noindent\proof  \emph{Step 1:  } Since $r_k\sin(2\pi(\theta_k-\gamma))$ is unbounded for every $\gamma\in \mathcal{I}_\Gamma$ and since $\mathbb{T}$ is compact we can, after passing to a subsequence, assume  that  $\{-z_k\}_{k\in\N}$ is contained in the interior of an escape cone $\mathcal{C}_{\alpha,\beta}$ for some $\beta\in \mathbb{T}\backslash \mathcal{I}_\Gamma$ and that  $\theta_k\rightarrow \theta^\ast\in[-\alpha+\beta,\alpha+\beta]$.

As the problem is rotation invariant, we may for simplicity assume that $\beta=0$. In the following we therefore omit the subscripts in $\eta_\beta, r_\beta$, etc.. As $\|z_k\|$ is unbounded, we can, after passing to yet another subsequence, write $z_k$ as 
$$
z_k= -\big(\eta (\theta^\ast)(n_k+v_k)+\rho(\theta^\ast)\big)\vec{\ell}( \theta_k),
$$  
where $n_k\in\N$ is  strictly increasing, $v_k\in [0,1]$ converges to $v^\ast$, and $\eta$, $\rho$ are the functions given by (A.iv). 
 If $z_k+\Gamma$ converges weakly to $\Gamma'$, then  $z_k+v_k\eta(\theta^\ast)\vec{\ell}(\theta_k)+\Gamma$  converges weakly to $v^\ast\eta(\theta^\ast)\vec{\ell}(\theta^\ast)+\Gamma'$.
As $z_k+v_k\eta(\theta^\ast)\vec{\ell}(\theta_k)$ still satisfies the  assumption of this lemma, we may further  simplify notation and assume
$$
z_k=-(\eta(\theta^\ast)k+\rho(\theta^\ast))\vec{\ell}(\theta_k).
$$
\emph{Step 2: }
We will later show that the following property is always satisfied: 
For every compact set $K\subset\R^2$ there exists $k^\ast\in\N$ such that  
$$
K\subset z_k+\mathcal{C}_{\alpha-|\theta_k|,\theta_k},\quad k\geq k^\ast.
$$
For $h\in C_c(\R^2)$, there exist $T_\theta,S_\theta>0$ and $t_\theta\in\R$ such that the parallelogram
$$ 
\widetilde{P}_{\theta}=\big\{(t_\theta +t)\vec{\ell}({\theta})+s\vec{d}(\theta): \ s\in[-S_\theta,S_\theta],\ t\in[0,T_\theta]\big\}
$$
satisfies $\text{supp}(h)\subset \widetilde{P}_{\theta}$.  Note that  $T_\theta,S_\theta$ can be chosen to depend continuously on $\theta\in[-\alpha,\alpha]$.
In order to simplify the argument, we define a slightly larger parallelogram $
P_{\theta}$ containing $
\widetilde{P}_{\theta}$.  To this end, let $m_{\theta},M_{\theta}\in\Z$ be the largest (respectively smallest) integer such that $ m_{\theta}\eta(\theta)\leq t_\theta$, and  $M_{\theta} \eta(\theta) \geq t_\theta +T_\theta$ and define
$$
P_{\theta}=\Big\{ (m_{\theta}-1/2+t)\eta(\theta) \vec{\ell}(\theta)+s\vec{d}(\theta): \ s\in[-S_\theta,S_\theta],\ t\in[0,M_{\theta}-m_{\theta}+1]\Big\},
$$
as well as the smaller parallelograms
$$
P_{\theta,m}=\Big\{ (m-1/2+t)\eta(\theta)  \vec{\ell}(\theta)+s\vec{d}(\theta): \ s\in[-S_\theta,S_\theta],\ t\in[0,1]\Big\},\quad m_{\theta}\leq m\leq M_{\theta}.
$$
Since all parameters defining $P_\theta$ are continuously depending on $\theta$, it first follows that the number of smaller parallelograms is bounded by $\max_{\theta\in[-\alpha,\alpha]}  T_\theta/\eta(\theta)+2$. Secondly, there exists a compact set $K$ such that $P_\theta\subset K$ for ever $\theta\in[-\alpha,\alpha]$. Consequently, by our previous assumption, for $k$ large enough we have that $ {P}_{\theta_k }\subset K \subset z_k+C_{\alpha-|\theta_k|,\theta_k }$.

\emph{Step 3:} For each $k,m$,  let $\psi_{k,m}:I_{k,m}\rightarrow\R^2$ be a re-parametrization by arc-length of the segment $\{z_k+\gamma(k+m+\theta)\}_{\theta\in[-\alpha,\alpha]}$ such that 
$$
\psi_{k,m}(0)=z_k+\gamma(k+m+\theta_k)=  \big(r(k+m+\theta_k) -k\eta(\theta^\ast)-\rho(\theta^\ast)\big)\vec{\ell}(\theta_k).
$$
For large values of $k$, $\psi_{k,m}(0)$ approximates $m\eta(\theta^\ast)\vec{\ell}(\theta^\ast)$.
Applying a first order Taylor approximation yields
\begin{equation}\label{eq:taylor}
\big\|\psi_{k,m}(t)-\big(r(k+m+\theta_k) -k\eta(\theta^\ast)-\rho(\theta^\ast)\big)\vec{\ell}(\theta_k)-t\psi'_{k,m}(0)\big\|\leq \frac{t^2}{2}\sup_{s\in I_{k,m}}\|\psi_{k,m}''(s)\|.
\end{equation}
Let us define $S_{\max}:= \max_{\theta\in[-\alpha,\alpha]}S_\theta$ (and $S_{\min}$ analogously). By conditions  (A.iii)-(A.v)  of Definition~\ref{def:spiraling} and
 $$\sup_{s\in I_{k,m}}\|\psi_{k,m}''(s)\|= \sup_{\theta\in [-\alpha,\alpha]} |\kappa(k+m+\theta)|\rightarrow 0, \quad\text{as }k\rightarrow\infty,$$
we have that for every $\delta>0$ there exists $k^\ast\in\N$ such that  \begin{enumerate}
\item[(i)] $|\kappa(k+m+\theta)|<\delta/8S_{\max}^2$, 
\item[(ii)]$
\big\|\psi^\prime_{k,m}(0)-\vec{d}(\theta_k)\big\|=\left\|\frac{\gamma^\prime(k+m+\theta_k)}{\|\gamma^\prime(k+m+\theta_k) \|}-\vec{d}(\theta_k)\right\|< {\delta}/{8S_{\max}},
$
\item[(iii)]
 $|r({k+m}+\theta^\ast)-\eta(\theta^\ast)(k+m)-\rho(\theta^\ast)|<\delta/4$, 
 \item[(iv)] $|r({k+m}+\theta_k)-r(k+m+\theta^\ast)|<\delta/4$,
 \end{enumerate}  whenever $k\geq k^\ast$. Note that we have used that the $C^2$-norm of $r$ is globally bounded to ensure the last estimate.
Now, if $t\in [-2 S_{\max},2S_{\max}]$ and $k\geq k^\ast$, then,  using \eqref{eq:taylor} and triangle inequality, we get
\begin{align*}
\big\|\psi_{k,m}(t)&-\eta (\theta^\ast)m\vec{\ell}(\theta_k)-t\vec{d}(\theta_k)\big\|
\\ 
\leq\ &|r(k+m+\theta^\ast)-r(k+m+\theta_k)| +|r(k+m+\theta^\ast)-\eta(\theta^\ast)(k+m)-\rho(\theta^\ast)|\times
\\
& +  |t|\|\psi'_{k,m}(0)-\vec{d}(\theta_k)\|+ \frac{t^2}{2}\sup_{s\in I_{k,m}}\|\psi_{k,m}''(s)\|
\\ 
<\ &\delta/4+\delta/4+\delta/4+\delta/4=\delta.
\end{align*}
Hence, if $\delta<\min\{S_{\min},\eta\}$, then $\psi_{k,m}\cap P_{\theta_k}\subset P_{\theta_k ,m}$, and $\psi_{k,m}(t) \notin P_{\theta_k,m}$ for every $|t|\geq 2S_{\max}$.

As $h$ is uniformly continuous we may choose $\delta$ according to $\varepsilon$ such that $\|x-y\|<\delta$ implies $|h(x)-h(y)|<\varepsilon$. Let $\lambda_k=\eta(\theta_k)\sin\big(\arccos(\vec{\ell}(\theta_k)\cdot\vec{d}(\theta_k))\big)$.
Then, for $k\geq k^\ast$
\begin{align*}
\left|\int_{z_k+\Gamma }h(x)\right.&\left.\ud\H^1(x)-\int_{ L_{\vec{d}(\theta_k),\lambda_k\Z}}h(x)\ud\H^1(x)\right| 
\\
&= \sum_{m=m_{1,\theta_k}}^{m_{2,\theta_k}}\left|\int_{(z_k+\Gamma)\cap P_{\theta_k,m}}h(x)\ud\H^1(x)-\int_{  L_{\vec{d}(\theta_k),\lambda_k\Z}\cap P_{\theta_k,m}}h(x)\ud\H^1(x)\right|
\\
&\leq\sum_{m=m_{1,\theta_k}}^{m_{2,\theta_k}}\int_{-2S_{\max}}^{2S_{\max}}\left|h(\psi_{k,m}(t)) -h\big( \eta(\theta_k)m\vec{\ell}({\theta_k})+\vec{d}(\theta_k)t\big)\right|\ud t
\\
&\leq 4 S_{\max} \sum_{m=m_{1,\theta_k}}^{m_{2,\theta_k}} \sup_{t\in[-2S_{\max},2S_{\max}]}\left|h(\psi_{k,m}(t)) -h\big( \eta(\theta_k)m\vec{\ell}({\theta_k})+\vec{d}(\theta_k)t\big)\right|
\\
& \leq 4S_{\max} (m_{2,\theta_k}-m_{1,\theta_k})\varepsilon\leq
4S_{\max} \left( \max_{\theta\in[-\alpha,\alpha]}  T_\theta/\eta(\theta)+2\right)\varepsilon.
\end{align*}
As $\theta_k\rightarrow\theta^\ast$  and consequently $\lambda_k\rightarrow\lambda^\ast$, it is clearly possible to find $k^\ast$ such that 
$$
\left|\int_{ L_{\vec{d}(\theta^\ast ),\lambda^\ast\Z}}h(x)\ud\H^1(x)-\int_{  L_{\vec{d}(\theta_k),\lambda_k\Z}}h(x)\ud\H^1(x)\right| \leq \varepsilon,
$$
whenever $k\geq k^\ast$. Therefore, by triangle inequality, the convergence $z_k+\Gamma\stackrel{w}{\rightarrow} L_{\vec{d}(\theta^\ast),\lambda^\ast\Z}$ follows.

\emph{Step 4:} It  remains to show that the assumption $K\subset z_k+\mathcal{C}_{\alpha-|\theta_k|,\theta_k}$  can always be satisfied for $k$ large enough. To this end, we distinguish the cases $|\theta^\ast|<\alpha$ and $|\theta^\ast|=\alpha$.
In the former case, we set $\alpha^\ast =(\alpha-|\theta^\ast|)/2 $ and observe that for $k$ large enough
 $-\frac{1}{2}\|z_k\|\vec{\ell}(\theta^\ast)+ \mathcal{C}_{\alpha^\ast,\theta^\ast}\subset z_k+ \mathcal{C}_{\alpha-|\theta_k|,\theta_k}$ and the left side of the inclusion relation will eventually cover any compactum.

If $|\theta^\ast|=\alpha$ and $\theta^\ast\notin\mathcal{I}_\Gamma$, then it is possible to apply a rotation with a small angle such that $\alpha$ is still the opening angle for the escape cone of $\beta=0$ and $|\theta^\ast|<\alpha$.

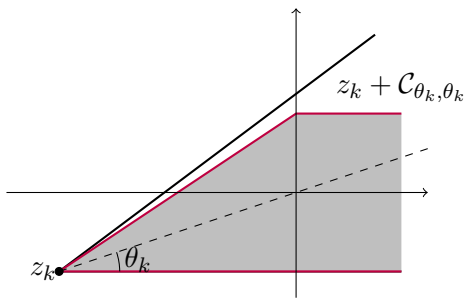
\begin{figure}[ht]
\begin{center}
\begin{tikzpicture}[scale=0.35]
\draw [thick] (-9,-3)--(4,-3);
\draw [thick] (-9,-3)--(3,6);
\path [fill=lightgray] (4,-3)--(-9,-3)--(0,3)--(4,3);
\draw [purple,thick] (4,-3)--(-9,-3)--(0,3)--(4,3);
\draw [->] (-11,0)--(5,0);
\draw [->] (0,-4)--(0,7);
\draw [dashed] (-9,-3)--(5,1.667);
\node at (-9.6,-3) {$z_k$};
\node at (4,4) {$z_k+\mathcal{C}_{\theta_k,\theta_k}$};
\draw [fill] (-9,-3) circle [radius=0.15];
\draw  (-6.7,-3) arc [radius=2.3, start angle=0, end angle= 18.5];
\node at (-6,-2.55) {$\theta_k$};
\end{tikzpicture}
\caption{The set from the left hand side of \eqref{eq:inclusion-cone} (gray shaded area) and the cone $z_k+\mathcal{C}_{\theta_k,\theta_k}$.} \label{fig:cone}
\end{center}
\end{figure}

If $|\theta^\ast|=\alpha$ and $\theta^\ast \in\mathcal{I}_\Gamma$,  then   one can rotate the problem so that  $\theta^\ast=0$. Moreover, we can without loss of generality assume that $\theta_k>0$. 
Let us consider the cone $\mathcal{C}_{\theta_k,\theta_k}$. A simple geometric argument  shows that 
\begin{equation}\label{eq:inclusion-cone}
\big\{(x,y)\in\R^2:\ x\geq -a_k,\ |y|\leq b_k,\text{ and } {2}b_k x+a_kb_k\geq a_ky \big\}\subset z_k+\mathcal{C}_{\theta_k,\theta_k},
\end{equation}
where $a_k=r_k\cos(2\pi \theta_k) $ and $b_k=r_k\sin(2\pi\theta_k)$, see Figure~\ref{fig:cone} for an illustration.
The assumptions that there is no unbounded subsequence of $\{z_k\}_k$ such that $r_{k_n}\sin(2\pi\theta_{k_n})$ is bounded  and $|\theta_k|\leq \alpha\leq 1/8$ then shows, after possibly passing to yet another subsequence, that   $ z_k+\mathcal{C}_{\theta_k,\theta_k}$ eventually covers any compactum. 

\emph{Step 5:} To see that every such collection of parallel lines is in fact a weak limit, we may assume that $\vec{d}=\lim_{\theta\searrow 0}\vec{d}_{\beta}(\theta)$ ($\vec{d}=\lim_{\theta\nearrow 0}\vec{d}_{\beta}(\theta)$ works exactly the same) and choose  
$$
z_k=\tau\vec{d}_\bot-(\eta_\beta(0)k+\rho_\beta(0))\vec{\ell}\big(\beta+1/\sqrt{k}\big),\ \beta\in\mathbb{T}\backslash\mathcal{I}_\Gamma.
$$ 
For this choice, the assumptions of the first part of this lemma are satisfied and one may repeat the arguments to show that $z_k+\Gamma\stackrel{w}{\rightarrow}\tau\vec{d}_\bot+L_{\vec{d},\eta\Z}$.\pbox

\begin{lemma}\label{lem:edges}
Let $ z_k =  -r_k \vec{\ell}(\theta_k)$, $k\in\N$, be an unbounded seqence. If there exists a pair $(\mathcal{N},\gamma),\ \mathcal{N}\subset \N$,  $\gamma\in\mathcal{I}_\Gamma$, such that 
$$
(i)\ \{z_k\}_{k \in\mathcal{N}}\text{ is unbounded\quad and \quad} (ii)\ \{r_{k}\sin(2\pi(\theta_{k}-\gamma))\}_{k \in \mathcal{N}}\text{ is  bounded,}
$$ then there exist a subsequence $\{z_{k_n}\}_{n\in\N}$,  and $z\in\R^2$,  such that 
$$
z_{k_n}+\Gamma\stackrel{w}{\rightarrow} z+E_{\gamma},
$$ 
where the set of \emph{parallel edges} $E_\gamma$ is given by
$$
E_{\gamma}:=\left\{\eta_\gamma k\vec{\ell}(\gamma)-t\vec{d}_{\gamma}^{-}:\ k\in\Z,\ t\in [0,\infty)\right\}\cup \left\{\eta_\gamma k\vec{\ell}(\gamma)+t\vec{d}_\gamma^+:\ k\in\Z,\ t\in [0,\infty)\right\},
$$
and $\eta_\gamma=\eta_\gamma(0)$, $\vec{d}_\gamma^-=\lim_{\theta\nearrow 0} \vec{d}_\gamma(\theta)$, and $\vec{d}_\gamma^+=\lim_{\theta\searrow 0} \vec{d}_\gamma(\theta)$.

Moreover, for every $\gamma\in\mathcal{I}_\Gamma$ and every $z\in\R^2$ there exist a sequence $\{z_k\}_{k\in\N}$ such that $z_k+\Gamma\stackrel{w}{\rightarrow}z+E_\gamma$.
\end{lemma}

\noindent\proof  After potentially passing to a subsequence, one can assume that $\mathcal{N}=\N$, and that $r_k$ is an increasing unbounded sequence. Moreover, by rotation invariance we may assume that $\gamma=0$. By assumption, the sequence $r_k\sin(2\pi\theta_k)$ is bounded which shows that there exists a  subsequence converging to $y^\ast$. Therefore, after passing to this subsequence, we see that if $z_k+\Gamma$  converges to $\Gamma'$, then $(-r_k\cos\theta_k,0)+\Gamma$ converges to $(0,-y^\ast)+\Gamma'$.
Therefore, we further simplify the problem and assume that $z_k=(-\widetilde{r}_k,0)$.
From here we can basically proceed as in the proof of Lemma~\ref{lem:lines1} with some minor adjustments.
 
Let us shortly point out where caution is needed.
For $h\in C_c(\R^2)$ the parallelograms need to be replaced by arrow shaped objects defined as follows: let $t^\ast\in\R$ and $S,T>0$ be chosen such that the   set
$$
\widetilde{A}_{h}:=\big\{(t^\ast+t)\vec{e}_1-s\vec{d}_0^-:\ t\in[0,T], s\in[0,S]\big\}
\ \cup \
\big\{(t^\ast+t)\vec{e}_1+s\vec{d}_0^+:\ t\in[0,T], s\in[0,S]\big\}
$$
contains $\text{supp}(h)$. Then one can proceed almost exactly as before by replacing $P_h$ by $A_h$, using the assumptions (B.iii)-(B.v) of Definition~\ref{def:spiraling} and applying two Taylor expansions for the right and left limits.

Again, we are left with showing that each such set of parallel edges is indeed a weak limit of translates. Setting $z_k=z-(\eta_\gamma(0)k+\rho_\gamma(0))\vec{\ell}(\gamma)$, $\gamma\in\mathcal{I}_\Gamma$, however yields that 
$z_k+\Gamma\stackrel{w}{\rightarrow} z+E_\gamma$. 
\pbox

\begin{theorem}\label{thm:weak-limits}
The set of weak limits of translates of a spiraling curve $\Gamma$ is given by 
\begin{equation}
\mathcal{W}_\Gamma=\mathcal{S}_\Gamma\cup   \mathcal{L}_\Gamma\cup\mathcal{E}_\Gamma,
\end{equation}
where $\mathcal{S}_\Gamma=\{z+\Gamma:\ z\in\R^2\}$, $\mathcal{E}_\Gamma= \big\{z+E_{\beta}:\ \beta\in  \mathcal{I}_\Gamma,\ z\in\R^2\big\}$, and 
$$
\mathcal{L}_\Gamma =\Big\{\tau\vec{d}_\bot +L_{\vec{d},\lambda\Z}:\ \tau\in\R,\ \vec{d}=\lim_{\theta\searrow 0}\vec{d}_{\beta}(\theta),\text{ or }\ \vec{d}=\lim_{\theta\nearrow 0}\vec{d}_{\beta}(\theta),\   \beta\in\mathbb{T}\Big\},
$$   
where $\lambda$ is defined as in Lemma~\ref{lem:lines1}.
\end{theorem}

\noindent \proof  It follows from Lemma~\ref{lem:lines1} and Lemma~\ref{lem:edges} that $\mathcal{S}_\Gamma\cup   \mathcal{L}_\Gamma\cup\mathcal{E}_\Gamma\subseteq \mathcal{W}_\Gamma$.

Now let $z_k+\Gamma\stackrel{w}{\rightarrow}\Gamma'$. If $\|z_{k}\|\leq C$, then there exist a converging subsequence $z_{k_n}\rightarrow z^\ast$ and it follows that $\Gamma'=z^\ast+\Gamma\in\mathcal{S}_\Gamma$. If $\{z_k\}_{k\in\N}$ is unbounded, then either there exist $\mathcal{N}\subset \N$ and $\gamma\in \mathcal{I}_\Gamma$  such that  $\{z_k\}_{k\in\mathcal{N}}$  is unbounded and $\{r_k \sin(2\pi(\theta_k-\gamma))\}_{k\in\mathcal{N}}$ is bounded, or not. Hence, the assumption of either Lemma~\ref{lem:lines1} or Lemma~\ref{lem:edges} are satisfied which leaves us with the limit $\Gamma'$ being either a set of parallel lines or a set of parallel edges, i.e.  $\mathcal{W}_\Gamma\subseteq \mathcal{S}_\Gamma\cup   \mathcal{L}_\Gamma\cup\mathcal{E}_\Gamma$. \pbox



\subsection{Spiraling curves as Gabor sampling trajectories }

After characterizing the set of all weak limits of translates we can now prove our main results. We first show that, for a certain class of windows,   spiraling curves are Gabor sampling trajectories if and only if they are uniqueness sets for the weighted modulation space $M^p_{1/\vartheta_s}(\R)$. 
Later, we then verify that certain spiraling curves are indeed such uniqueness sets.

\begin{theorem}\label{thm:sampling=uniqueness}
Let $\Gamma\in\R^2$ be a spiraling curve, $g\in \text{span}\{h_n:   n\in\N_0\}$, and $1\leq p<\infty$. Then $\Gamma$ is a Gabor sampling trajectory for $M^p(\R)$ if and only if $\Gamma$ is a uniqueness set for $M^p_{1/\vartheta_s}(\R)$ for some $s>2$.
\end{theorem}

\noindent To show this theorem we first need to state two auxiliary lemmas. The first one is a consequence of \cite[Corollary~3.9(c) \& Proposition~2.2]{grzi04}.
\begin{lemma}\label{lem:aux-hermite-coeff}
 If there exist $\alpha>0$ such that  $|g(t)|\lesssim e^{-\alpha|t|}$ and $|\widehat{g}(\xi)|\lesssim   e^{- \alpha  |\xi|}$, then   
$
|V_gg(z)|\lesssim  {1}/{\vartheta_{s}(z)}$, for every $s>0$.

\end{lemma}
The following lemma is a simple consequence of \cite[Corollary 26]{Ascensi}.

\begin{lemma} \label{lem:ascensi-analytic}
Let $g\in \text{span}\{h_n:   n\in\N_0\}$.   For every $f\in M^p_{1/\vartheta_s}(\R)$, $1\leq p<\infty$, $z\in\R^2$, and $\beta\in\mathbb{T}$, it holds that $V_{g}f(z+t\vec{\ell}(\beta))$ is a real analytic function in   $t\in\R$.
\end{lemma}
\proof 
First, observe that  $V_{g}f(z+t\vec{\ell}(\beta)) $ is real analytic if and only if 
$$
V_{g}f\big(R( \beta )R( -\beta)(z+t\vec{\ell}(\beta))\big)=V_{\mu( \beta )g}(\mu(\beta)f)\big(R(-\beta)z+t\vec{e}_1\big)
$$ 
is real analytic.  For $\beta\in\mathbb{T}$ one has that $\mu(\beta)g\in \mathcal{M}_s^1(\R)$ by Lemma~\ref{lem:aux-hermite-coeff}  and $|\mu(\beta)g(t)|\lesssim e^{-\alpha|t|}$. Since the metaplectic rotation leaves $M^p_{1/\vartheta_s}(\R)$  invariant,   it follows  by \cite[Corollary 26]{Ascensi}  that $V_{\mu( \beta)g}(\mu(\beta)f)\big(R(-\beta)z+t\vec{e}_1\big)$ is real analytic. \pbox \medskip\medskip

\noindent \textbf{Proof of Theorem~\ref{thm:sampling=uniqueness}:} First, observe that  $g\in\mathcal{M}^1_s(\R)$ by Lemma~\ref{lem:aux-hermite-coeff}.  We may therefore apply Theorem~\ref{thm:ascensi}. By Theorem~\ref{thm:weak-limits}, the set of weak limits of translates consists of finite shifts of $\Gamma$ and the sets of parallel lines and edges. As $V_gf(z-w)=e^{2\pi i \xi y}V_g\pi(w)f(z)$, and $M_{1/\vartheta_s}^p(\R)$ is invariant under time-frequency shifts, it follows that $z+\Gamma$ is a uniqueness set for $M_{1/\vartheta_s}^p(\R)$ if and only if $\Gamma$ is a uniqueness set for $M_{1/\vartheta_s}^p(\R)$.

We now show that every set of parallel lines (resp. parallel edges) is a  uniqueness set. To show the two cases at once, we prove that any collection of parallel half lines $\{\eta k\vec{d}_1 +t\vec{d}_2:\ t\geq 0\}$, $\vec{d_1},\vec{d}_2$ not colinear, and $\eta>0$, is  a uniqueness set for $M_{1/\vartheta_s}^p(\R)$.   
Since $V_gf$ restricted to a line in $\R^2$ is real analytic by  Lemma~\ref{lem:ascensi-analytic}, it follows that  $V_g f|_{\{\eta k\vec{d}_1 +t\vec{d}_2:\ t\geq 0, \ k\in\Z\}}=0$ implies $V_g f|_{\{\eta k\vec{d}_1 +t\vec{d}_2:\ t\in\R, \ k\in\Z\}}=0$. 

The function  $\mu(\theta)g$ is also a  linear combination of Hermite functions and thus  has  only finitely many zeros. Therefore, it follows that $\bigcup_{k\in\Z}\widetilde{\text{supp}}(T_{\delta k}\mu(\theta)g)=\R$ for every $\delta>0$.
Finally, as $g\in A_{\vartheta_s}(\R)$, we may apply  Proposition~\ref{prop:lines-uniqueness} to  see that $V_g f|_{\{\eta k\vec{d}_1 +t\vec{d}_2:\ t\in\R, \ k\in\Z\}}=0$ implies $f=0$.
\pbox

\begin{corollary}
Let $g\in \text{span}\{h_n:   n\in\N_0\}$,  $1\leq p<\infty$, $\eta>0$. Moreover, let    $P\subset \R^2$   and $\{z_1,...,z_n\}\subset \R^2$  be such that the assumptions $(a)-(b)$ in Theorem~\ref{thm:main} are satisfied.  The  collection of star shaped polygons $P_\eta$ and $\mathcal{S}(z_1,...,z_n)$  are Gabor sampling trajectory for $M^p(\R)$.
\end{corollary}
\proof Take a line $\ell\subset \R^2$ that contains one of the edges of $P$ (resp. that contains the line segment $s(z_1,z_2)$) and consider the collection $\{\eta k\ell\}_{k\in\N}$. If $V_gf|_{P_\eta}=0$ (resp. $V_gf|_{\mathcal{S}(z_1,...,z_n)}=0$), then, as $V_gf$ is real analytic on any  line $\eta k\ell$ and zero on a subset of nonzero measure, it follows that $V_gf|_{\{\eta k\ell\}_{k\in\N}}=0$. Arguing as before, we get
$\bigcup_{k\in\N} \widetilde{\text{supp}}(T_{k\eta}\mu(\beta)g)=\R$. Therefore, it follows that
that $f=0$. \pbox

\begin{corollary}
Let  $g\in \text{span}\{h_n:   n\in\N_0\}$, $\eta>0$, and $1\leq p<\infty$. The collection of concentric circles $O_\eta$ is a Gabor sampling trajectory for $M^p(\R)$.
\end{corollary} 
\proof 
For $g=\sum_{k=0}^n\alpha_k h_n$,   we denote by $F$ the polyanalytic function of order $n$ given by  $F:=\sum_{k=0}^n\alpha_k F_k$, where    the functions $F_0,...,F_n$ are given in \eqref{def-poly-bargmann-trafo}. Then $V_g f(\overline{z})=0$ exactly when $F(z)=0$.  Multiplying $F$ by $z^n$  changes the zero set of $F$ at most in the origin, and results in  a reduced polyanalytic function. By assumption, and since $\overline{O_\eta}=O_\eta$ we hence have that $z^nF|_{O_\eta}=0$.
By Lemma~\ref{eq:cauchy-formula}, it thus follows that $F(z)=0$ for every $z\in B_{R_0}(0)$, where $R_0$ can be chosen arbitrarily large. Consequently, $F(z)=0$ for every $z\in\C$ which implies that $f=0$  and $O_\eta$ is a uniqueness set for $M^p_{1/\vartheta_s}(\R)$. \pbox

\section*{Acknowledgement}
The author would like to thank Felipe Negrera and Jos{\'e}-Luis Romero for valuable discussions and suggesting some of the references. The author also acknowledges the support of the Austrian Science Fund (FWF) through the Erwin Schr{\"o}dinger Fellowship J-4254.

\bibliographystyle{plain}
\bibliography{paperbib}

\end{document}